\documentclass{article}
\usepackage{amsmath}
\usepackage{amssymb}
\usepackage{amsfonts}
\usepackage{hyperref}

\newtheorem{theorem}{Theorem}

\newtheorem{proposition}[theorem]{Proposition}
\newtheorem{remark}[theorem]{Remark}

\newenvironment{proof}[1][Proof]{\noindent\textbf{#1.} }{\ \rule{0.5em}{0.5em}}
\input{tcilatex}
\numberwithin{theorem}{section}
\numberwithin{equation}{section}

\begin{document}

\title{$1$-Jet Riemann--Finsler Geometry for the Three-Dimensional Time}
\author{Gheorghe Atanasiu and Mircea Neagu \\
{\scriptsize February 2010; Revised September 2010}\\
{\scriptsize (correction upon the canonical nonlinear connection)}}
\date{}
\maketitle

\begin{abstract}
The aim of this paper is to develop on the $1$-jet space $J^{1}(\mathbb{R}%
,M^{3})$ the Finsler-like geometry (in the sense of distinguished (d-)
connection, d-torsions and d-curvatures) of the rheonomic Berwald-Mo\'{o}r
metric of order three 
\begin{equation*}
\mathring{F}(t,y)=\sqrt{h^{11}(t)}\cdot \sqrt[3]{y_{1}^{1}y_{1}^{2}y_{1}^{3}}%
.
\end{equation*}%
Some natural geometrical field theories (gravitational and electromagnetic)
produced by the above rheonomic Berwald-Mo\'{o}r metric are also exposed.
\end{abstract}



\section{Introduction}

\hspace{5mm}It is an well known fact that, in order to create the Relativity
Theory, Einstein was forced to use the Riemannian geometry instead of the
classical Euclidean geometry, the first one representing the natural
ma\-the\-ma\-ti\-cal model for the local \textit{isotropic} space-time. But,
there are recent studies of physicists which suggest a \textit{non-isotropic}
perspective of the space-time. For example, in Pavlov's opinion \cite%
{Pavlov-2}, the concept of inertial body mass emphasizes the necessity of
study of local non-isotropic spaces. Obviously, for the study of
non-isotropic physical phenomena, the Finsler geometry is very useful as
ma\-the\-ma\-ti\-cal framework.

The studies of Russian scholars (Asanov \cite{Asanov[1]}, Garas'ko \cite%
{Garasko-Book} and Pavlov \cite{Garasko-Pavlov}, \cite{Pavlov-1}, \cite%
{Pavlov-2}) emphasize the importance of the Finsler geometry which is
characterized by the total equality of all non-isotropic directions. For
such a reason, Asanov, Pavlov and their co-workers underline the important
role played in the theory of space-time structure and gravitation, as well
as in unified gauge field theories, by the Berwald-Mo\'{o}r metric (whose
certain Finsler geometrical properties are studied by Matsumoto and Shimada
in the works \cite{Matsumoto}, \cite{Mats-Shimada}, \cite{Shimada})%
\begin{equation*}
F:TM\rightarrow \mathbb{R},\mathbb{\qquad }F(y)=\left(
y^{1}y^{2}...y^{n}\right) ^{1/n}.
\end{equation*}%
Because any of such directions can be related to the proper time of an
inertial reference frame, Pavlov considers that it is appropriate as such
spaces to be generically called \textit{"multi-dimensional time"} \cite%
{Pavlov-2}. In the framework of the $3$- and $4$-dimensional linear space
with Berwald-Mo\'{o}r metric (i.e. the three- and four-dimensional time),
Pavlov and his co-workers \cite{Garasko-Pavlov}, \cite{Pavlov-1}, \cite%
{Pavlov-2} offer some new physical approaches and geometrical
interpretations such as:

1. physical events = points in the multi-dimensional time;

2. straight lines = shortest curves;

3. intervals = distances between the points along of a straight line;

4. light pyramids $\Leftrightarrow $ light cones in a pseudo-Euclidian space;

5. simultaneous surfaces = the surfaces of simultaneous physical events.

According to Olver's opinion \cite{Olver}, we appreciate that the $1$-jet
fibre bundle is a basic object in the study of classical and quantum field
theories. For such geometrical and physical reasons, this paper is devoted
to the development on the $1$-jet space $J^{1}(\mathbb{R},M^{3})$ of the
Finsler-like geometry (together with a theoretical-geometric gravitational
and electromagnetic field theory) of the rheonomic Berwald-Mo\'{o}r metric%
\textit{\ }%
\begin{equation*}
\mathring{F}:J^{1}(\mathbb{R},M^{3})\rightarrow \mathbb{R},\qquad \mathring{F%
}(t,y)=\sqrt{h^{11}(t)}\cdot \sqrt[3]{y_{1}^{1}y_{1}^{2}y_{1}^{3}},
\end{equation*}%
where $h_{11}(t)$ is a Riemannian metric on $\mathbb{R}$ and $%
(t,x^{1},x^{2},x^{3},y_{1}^{1},y_{1}^{2},y_{1}^{3})$ are the coordinates of
the $1$-jet space $J^{1}(\mathbb{R},M^{3})$.

The geometry (in the sense of d-connections, d-torsions, d-curvatures,
gravitational and electromagnetic geometrical theories) produced by an
arbitrary jet rheonomic Lagrangian function $L:J^{1}(\mathbb{R}%
,M^{n})\rightarrow \mathbb{R}$ is now completely done in the second author's
paper \cite{Neagu-Rheon}. We point out that the geometrical ideas from \cite%
{Neagu-Rheon} are similar, but however distinct ones, with those exposed by
Miron and Anastasiei in the classical Lagrangian geometry \cite{Mir-An}. In
fact, the geometrical ideas from \cite{Neagu-Rheon} (which we called the jet
geometrical theory of the \textit{relativistic rheonomic Lagrange spaces})
were initially stated by Asanov in \cite{Asanov[2]} and developed further by
the second author of this paper in the book \cite{Neagu Carte}.

In the sequel, we apply the general geometrical results from \cite%
{Neagu-Rheon} to the particular rheonomic Berwald-Mo\'{o}r metric\textit{\ }$%
\mathring{F}$, in order to obtain what we called the $1$\textit{-jet
Riemann-Finsler geometry of the three-dimensional time}.

\section{Preliminary notations and formulas}

\hspace{5mm}Let $(\mathbb{R},h_{11}(t))$ be a Riemannian manifold, where $%
\mathbb{R}$ is the set of real numbers. The Christoffel symbol of the
Riemannian metric $h_{11}(t)$ is%
\begin{equation*}
\varkappa _{11}^{1}=\frac{h^{11}}{2}\frac{dh_{11}}{dt},\qquad h^{11}=\frac{1%
}{h_{11}}>0.
\end{equation*}%
Let also $M^{3}$ be a manifold of dimension three, whose local coordinates
are $(x^{1},x^{2},x^{3})$. Let us consider the $1$-jet space $J^{1}(\mathbb{R%
},M^{3})$, whose local coordinates are%
\begin{equation*}
(t,x^{1},x^{2},x^{3},y_{1}^{1},y_{1}^{2},y_{1}^{3}).
\end{equation*}%
These transform by the rules (the Einstein convention of summation is used
throughout this work):%
\begin{equation}
\widetilde{t}=\widetilde{t}(t),\quad \widetilde{x}^{p}=\widetilde{x}%
^{p}(x^{q}),\quad \widetilde{y}_{1}^{p}=\dfrac{\partial \widetilde{x}^{p}}{%
\partial x^{q}}\dfrac{dt}{d\widetilde{t}}\cdot y_{1}^{q},\qquad p,q=%
\overline{1,3},  \label{tr-rules}
\end{equation}%
where $d\widetilde{t}/dt\neq 0$ and rank $(\partial \widetilde{x}%
^{p}/\partial x^{q})=3$. We consider that the manifold $M^{3}$ is endowed
with a tensor of kind $(0,3)$, given by the local components $G_{pqr}(x)$,
which is totally symmetric in the indices $p$, $q$ and $r$. Suppose that the
d-tensor 
\begin{equation*}
G_{ij1}=6G_{ijp}y_{1}^{p},
\end{equation*}%
is non-degenerate, that is there exists the d-tensor $G^{jk1}$ on $J^{1}(%
\mathbb{R},M^{3})$ such that $G_{ij1}G^{jk1}=\delta _{i}^{k}.$

In this geometrical context, if we use the notation $%
G_{111}=G_{pqr}y_{1}^{p}y_{1}^{q}y_{1}^{r}$, we can consider the \textit{%
third-root Finsler-like function} \cite{Shimada}, \cite{Atanasiu-Neagu} (it
is $1$-positive homogenous in the variable $y$):%
\begin{equation}
F(t,x,y)=\sqrt[3]{G_{pqr}(x)y_{1}^{p}y_{1}^{q}y_{1}^{r}}\cdot \sqrt{h^{11}(t)%
}=\sqrt[3]{G_{111}(x,y)}\cdot \sqrt{h^{11}(t)},  \label{F}
\end{equation}%
where the Finsler function $F$ has as domain of definition all values $%
(t,x,y)$ which verify the condition $G_{111}(x,y)\neq 0.$ If we denote $%
G_{i11}=3G_{ipq}y_{1}^{p}y_{1}^{q}$, then the $3$-positive homogeneity of
the "$y$-function" $G_{111}$ (this is in fact a d-tensor on the $1$-jet
space $J^{1}(\mathbb{R},M^{3})$) leads to the equalities:%
\begin{equation*}
G_{i11}=\frac{\partial G_{111}}{\partial y_{1}^{i}},\quad
G_{i11}y_{1}^{i}=3G_{111},\quad G_{ij1}y_{1}^{j}=2G_{i11},
\end{equation*}%
\begin{equation*}
G_{ij1}=\frac{\partial G_{i11}}{\partial y_{1}^{j}}=\frac{\partial
^{2}G_{111}}{\partial y_{1}^{i}\partial y_{1}^{j}},\quad
G_{ij11}y_{1}^{i}y_{1}^{j}=6G_{111},\quad \frac{\partial G_{ij1}}{\partial
y_{1}^{k}}=6G_{ijk}.
\end{equation*}

The \textit{fundamental metrical d-tensor} produced by $F$ is given by the
formula%
\begin{equation*}
g_{ij}(t,x,y)=\frac{h_{11}(t)}{2}\frac{\partial ^{2}F^{2}}{\partial
y_{1}^{i}\partial y_{1}^{j}}.
\end{equation*}%
By direct computations, the fundamental metrical d-tensor takes the form%
\begin{equation}
g_{ij}(x,y)=\frac{G_{111}^{-1/3}}{3}\left[ G_{ij1}-\frac{1}{3G_{111}}%
G_{i11}G_{j11}\right] .  \label{g-(ij)-general}
\end{equation}%
Moreover, taking into account that the d-tensor $G_{ij1}$ is non-degenerate,
we deduce that the matrix $g=(g_{ij})$ admits the inverse $g^{-1}=(g^{jk})$.
The entries of the inverse matrix $g^{-1}$ are%
\begin{equation}
g^{jk}=3G_{111}^{1/3}\left[ G^{jk1}+\frac{G_{1}^{j}G_{1}^{k}}{3\left(
G_{111}-\mathcal{G}_{111}\right) }\right] ,  \label{g+(jk)-general}
\end{equation}%
where $G_{1}^{j}=G^{jp1}G_{p11}$ and $3\mathcal{G}%
_{111}=G^{pq1}G_{p11}G_{q11}.$

\section{The rheonomic Berwald-Mo\'{o}r metric}

\hspace{5mm}Beginning with this Section we will focus only on the \textit{%
rheonomic Berwald-Mo\'{o}r metric of order three}, which is the Finsler-like
metric (\ref{F}) for the particular case%
\begin{equation*}
G_{pqr}=\left\{ 
\begin{array}{ll}
\dfrac{1}{3!}, & \{p,q,r\}\text{ - distinct indices}\medskip \\ 
0, & \text{otherwise.}%
\end{array}%
\right.
\end{equation*}%
Consequently, the rheonomic Berwald-Mo\'{o}r metric of order three is given
by%
\begin{equation}
\mathring{F}(t,y)=\sqrt{h^{11}(t)}\cdot \sqrt[3]{y_{1}^{1}y_{1}^{2}y_{1}^{3}}%
.  \label{rheon-B-M}
\end{equation}%
Moreover, using preceding notations and formulas, we obtain the following
relations:%
\begin{equation*}
G_{111}=y_{1}^{1}y_{1}^{2}y_{1}^{3},\quad G_{i11}=\frac{G_{111}}{y_{1}^{i}},
\end{equation*}%
\begin{equation*}
G_{ij1}=\left( 1-\delta _{ij}\right) \frac{G_{111}}{y_{1}^{i}y_{1}^{j}}\text{
(no sum by }i\text{ or }j\text{),}
\end{equation*}%
where $\delta _{ij}$ is the Kronecker symbol. Because we have%
\begin{equation*}
\det \left( G_{ij1}\right) _{i,j=\overline{1,3}}=2G_{111}\neq 0,
\end{equation*}%
we find%
\begin{equation*}
G^{jk1}=\frac{(1-2\delta ^{jk})}{2G_{111}}y_{1}^{j}y_{1}^{k}\text{ (no sum
by }j\text{ or }k\text{).}
\end{equation*}%
It follows that we have $\mathcal{G}_{111}=(1/2)G_{111}$ and $%
G_{1}^{j}=(1/2)y_{1}^{j}$.

Replacing now the preceding computed entities into the formulas (\ref%
{g-(ij)-general}) and (\ref{g+(jk)-general}), we get%
\begin{equation}
g_{ij}=\frac{\left( 2-3\delta _{ij}\right) }{9}\frac{G_{111}^{2/3}}{%
y_{1}^{i}y_{1}^{j}}\text{ (no sum by }i\text{ or }j\text{)}
\label{g-jos-(ij)}
\end{equation}%
and%
\begin{equation}
g^{jk}=(2-3\delta ^{jk})G_{111}^{-2/3}y_{1}^{j}y_{1}^{k}\text{ (no sum by }j%
\text{ or }k\text{).}  \label{g-sus-(jk)}
\end{equation}

Using a general formula from the paper \cite{Neagu-Rheon}, we find the
following geometrical result:

\begin{theorem}
For the rheonomic Berwald-Mo\'{o}r metric (\ref{rheon-B-M}), the \textit{%
energy action functional}%
\begin{equation*}
\mathbb{\mathring{E}}(t,x(t))=\int_{a}^{b}\sqrt[3]{\left\{
y_{1}^{1}y_{1}^{2}y_{1}^{3}\right\} ^{2}}\cdot h^{11}\sqrt{h_{11}}dt
\end{equation*}%
produces on the $1$-jet space $J^{1}(\mathbb{R},M^{3})$ the \textit{%
canonical nonlinear connection}%
\begin{equation}
\Gamma =\left( M_{(1)1}^{(i)}=-\varkappa _{11}^{1}y_{1}^{i},\text{ }%
N_{(1)j}^{(i)}=0\right) .  \label{can-nlc=0}
\end{equation}
\end{theorem}

Because the canonical nonlinear connection (\ref{can-nlc=0}) has the spatial
components equal to zero, it follows that our subsequent geometrical theory
becomes trivial, in a way. For such a reason, in order to avoid the
triviality of our theory and in order to have a certain kind of symmetry, we
will use on the $1$-jet space $J^{1}(\mathbb{R},M^{3})$, by an "a priori"
definition, the following nonlinear connection:%
\begin{equation}
\mathring{\Gamma}=\left( M_{(1)1}^{(i)}=-\varkappa _{11}^{1}y_{1}^{i},\text{ 
}N_{(1)j}^{(i)}=-\frac{\varkappa _{11}^{1}}{2}\delta _{j}^{i}\right) .
\label{nlc-B-M}
\end{equation}

\section{Cartan canonical connection. d-Torsions and d-curvatures}

\hspace{5mm}The importance of the nonlinear connection (\ref{nlc-B-M}) is
coming from the possibility of construction of the dual \textit{adapted bases%
} of d-vector fields%
\begin{equation}
\left\{ \frac{\delta }{\delta t}=\frac{\partial }{\partial t}+\varkappa
_{11}^{1}y_{1}^{p}\frac{\partial }{\partial y_{1}^{p}},\text{ }\frac{\delta 
}{\delta x^{i}}=\frac{\partial }{\partial x^{i}}+\frac{\varkappa _{11}^{1}}{2%
}\frac{\partial }{\partial y_{1}^{i}},\text{ }\dfrac{\partial }{\partial
y_{1}^{i}}\right\} \subset \mathcal{X}(E)  \label{a-b-v}
\end{equation}%
and d-covector fields%
\begin{equation}
\left\{ dt,\text{ }dx^{i},\text{ }\delta y_{1}^{i}=dy_{1}^{i}-\varkappa
_{11}^{1}y_{1}^{i}dt-\frac{\varkappa _{11}^{1}}{2}dx^{i}\right\} \subset 
\mathcal{X}^{\ast }(E),  \label{a-b-co}
\end{equation}%
where $E=J^{1}(\mathbb{R},M^{3})$. Note that, under a change of coordinates (%
\ref{tr-rules}), the elements of the adapted bases (\ref{a-b-v}) and (\ref%
{a-b-co}) transform as classical tensors. Consequently, all subsequent
geometrical objects on the $1$-jet space $J^{1}(\mathbb{R},M^{3})$ (such as
Cartan canonical connection, torsion, curvature etc.) will be described in
local adapted components.

Using a general result from \cite{Neagu-Rheon}, by direct computations, we
can give the following important geometrical result:

\begin{theorem}
The Cartan canonical $\mathring{\Gamma}$-linear connection, produced by the
rheonomic Berwald-Mo\'{o}r metric (\ref{rheon-B-M}), has the following
adapted local components:%
\begin{equation*}
C\mathring{\Gamma}=\left( \varkappa _{11}^{1},\text{ }G_{j1}^{k}=0,\text{ }%
L_{jk}^{i}=\frac{\varkappa _{11}^{1}}{2}C_{j(k)}^{i(1)},\text{ }%
C_{j(k)}^{i(1)}\right) ,
\end{equation*}%
where, if we use the notation%
\begin{equation*}
A_{jk}^{i}=\frac{3\delta _{j}^{i}+3\delta _{k}^{i}+3\delta _{jk}-9\delta
_{j}^{i}\delta _{jk}-2}{9}\text{ (no sum by }i,\text{ }j\text{ or }k\text{),}
\end{equation*}%
then%
\begin{equation*}
C_{j(k)}^{i(1)}=A_{jk}^{i}\cdot \frac{y_{1}^{i}}{y_{1}^{j}y_{1}^{k}}\text{
(no sum by }i,\text{ }j\text{ or }k\text{).}
\end{equation*}
\end{theorem}

\begin{proof}
Via the Berwald-Mo\'{o}r derivative operators (\ref{a-b-v}) and (\ref{a-b-co}%
), we use the general formulas which give the adapted components of the
Cartan canonical connection, namely \cite{Neagu-Rheon}%
\begin{equation*}
G_{j1}^{k}=\frac{g^{km}}{2}\frac{\delta g_{mj}}{\delta t},\quad L_{jk}^{i}=%
\frac{g^{im}}{2}\left( \frac{\delta g_{jm}}{\delta x^{k}}+\frac{\delta g_{km}%
}{\delta x^{j}}-\frac{\delta g_{jk}}{\delta x^{m}}\right) ,
\end{equation*}%
\begin{equation*}
C_{j(k)}^{i(1)}=\frac{g^{im}}{2}\left( \frac{\partial g_{jm}}{\partial
y_{1}^{k}}+\frac{\partial g_{km}}{\partial y_{1}^{j}}-\frac{\partial g_{jk}}{%
\partial y_{1}^{m}}\right) =\frac{g^{im}}{2}\frac{\partial g_{jk}}{\partial
y_{1}^{m}}.
\end{equation*}
\end{proof}

\begin{remark}
The below properties of the d-tensor $C_{j(k)}^{i(1)}$ are true (sum by $m$):%
\begin{equation}
C_{j(k)}^{i(1)}=C_{k(j)}^{i(1)},\quad C_{j(m)}^{i(1)}y_{1}^{m}=0,\quad
C_{j(m)}^{m(1)}=0.  \label{equalitie-C}
\end{equation}%
For similar properties, please see also the papers \cite{Atanasiu-Neagu}, 
\cite{Mats-Shimada}, \cite{Neagu-B-M-4} or \cite{Shimada}.
\end{remark}

\begin{remark}
The coefficients $A_{ij}^{l}$ have the following values:%
\begin{equation}
A_{ij}^{l}=\left\{ 
\begin{array}{ll}
-\dfrac{2}{9}, & i\neq j\neq l\neq i\medskip \\ 
\dfrac{1}{9}, & i=j\neq l\text{ or }i=l\neq j\text{ or }j=l\neq i\medskip \\ 
-\dfrac{2}{9}, & i=j=l.%
\end{array}%
\right.  \label{A-(ijk)}
\end{equation}
\end{remark}

\begin{theorem}
The Cartan canonical connection $C\mathring{\Gamma}$ of the rheonomic
Berwald-Mo\'{o}r metric (\ref{rheon-B-M}) has \textbf{three} effective
adapted local torsion d-tensors:%
\begin{equation*}
\begin{array}{c}
P_{(1)i(j)}^{(k)\text{ }(1)}=-\dfrac{\varkappa _{11}^{1}}{2}%
C_{i(j)}^{k(1)},\quad P_{i(j)}^{k(1)}=C_{i(j)}^{k(1)},\medskip \\ 
R_{(1)1j}^{(k)}=\dfrac{1}{2}\left[ \dfrac{d\varkappa _{11}^{1}}{dt}%
-\varkappa _{11}^{1}\varkappa _{11}^{1}\right] \delta _{j}^{k}.%
\end{array}%
\end{equation*}
\end{theorem}

\begin{proof}
A general $h$-normal $\Gamma $-linear connection on the 1-jet space $J^{1}(%
\mathbb{R},M^{3})$ is characterized by \textit{eight} effective d-tensors of
torsion (for more details, please see \cite{Neagu-Rheon}). For our Cartan
canonical connection $C\mathring{\Gamma}$ these reduce to the following 
\textit{three} (the other five cancel):%
\begin{equation*}
{P_{(1)i(j)}^{(k)\text{ }(1)}={\dfrac{\partial N_{(1)i}^{(k)}}{\partial
y_{1}^{j}}}-L_{ji}^{k}},\quad R_{(1)1j}^{(k)}={\dfrac{\delta M_{(1)1}^{(k)}}{%
\delta x^{j}}}-{\dfrac{\delta N_{(1)j}^{(k)}}{\delta t}},\quad
P_{i(j)}^{k(1)}=C_{i(j)}^{k(1)}.
\end{equation*}
\end{proof}

\begin{theorem}
The Cartan canonical connection $C\mathring{\Gamma}$ of the rheonomic
Berwald-Mo\'{o}r metric (\ref{rheon-B-M}) has \textbf{three} effective
adapted local curvature d-tensors:%
\begin{equation*}
\begin{array}{c}
R_{ijk}^{l}=\dfrac{\varkappa _{11}^{1}\varkappa _{11}^{1}}{4}%
S_{i(j)(k)}^{l(1)(1)},\quad P_{ij(k)}^{l\text{ }(1)}=\dfrac{\varkappa
_{11}^{1}}{2}S_{i(j)(k)}^{l(1)(1)},\medskip \\ 
S_{i(j)(k)}^{l(1)(1)}={{\dfrac{\partial C_{i(j)}^{l(1)}}{\partial y_{1}^{k}}}%
-{\dfrac{\partial C_{i(k)}^{l(1)}}{\partial y_{1}^{j}}}%
+C_{i(j)}^{m(1)}C_{m(k)}^{l(1)}-C_{i(k)}^{m(1)}C_{m(j)}^{l(1)}.}%
\end{array}%
\end{equation*}
\end{theorem}

\begin{proof}
A general $h$-normal $\Gamma $-linear connection on the 1-jet space $J^{1}(%
\mathbb{R},M^{3})$ is characterized by \textit{five} effective d-tensors of
curvature (for more details, please see \cite{Neagu-Rheon}). For our Cartan
canonical connection $C\mathring{\Gamma}$ these reduce to the following 
\textit{three} (the other two cancel):%
\begin{equation*}
{R_{ijk}^{l}={\dfrac{\delta L_{ij}^{l}}{\delta x^{k}}}-{\dfrac{\delta
L_{ik}^{l}}{\delta x^{j}}}+L_{ij}^{m}L_{mk}^{l}-L_{ik}^{m}L_{mj}^{l},}
\end{equation*}%
\begin{equation*}
{P_{ij(k)}^{l\;\;(1)}={\dfrac{\partial L_{ij}^{l}}{\partial y_{1}^{k}}}%
-C_{i(k)|j}^{l(1)}+C_{i(m)}^{l(1)}P_{(1)j(k)}^{(m)\;\;(1)},}
\end{equation*}%
\begin{equation*}
S_{i(j)(k)}^{l(1)(1)}={{\dfrac{\partial C_{i(j)}^{l(1)}}{\partial y_{1}^{k}}}%
-{\dfrac{\partial C_{i(k)}^{l(1)}}{\partial y_{1}^{j}}}%
+C_{i(j)}^{m(1)}C_{m(k)}^{l(1)}-C_{i(k)}^{m(1)}C_{m(j)}^{l(1)},}
\end{equation*}%
where%
\begin{equation*}
{C_{i(k)|j}^{l(1)}=}\frac{\delta {C_{i(k)}^{l(1)}}}{\delta x^{j}}+{%
C_{i(k)}^{m(1)}L_{mj}^{l}}-{C_{m(k)}^{l(1)}L_{ij}^{m}}-{%
C_{i(m)}^{l(1)}L_{kj}^{m}}.
\end{equation*}
\end{proof}

\begin{remark}
The curvature d-tensor $S_{i(j)(k)}^{l(1)(1)}$ has the properties%
\begin{equation*}
S_{i(j)(k)}^{l(1)(1)}+S_{i(k)(j)}^{l(1)(1)}=0,\quad S_{i(j)(j)}^{l(1)(1)}=0%
\text{ (no sum by }j\text{).}
\end{equation*}
\end{remark}

\begin{theorem}
The expressions of the curvature d-tensor $S_{i(j)(k)}^{l(1)(1)}$ are:

\begin{enumerate}
\item $S_{i(i)(k)}^{l(1)(1)}=-\dfrac{1}{9}\dfrac{y_{1}^{l}}{\left(
y_{1}^{i}\right) ^{2}y_{1}^{k}}$ ($i\neq k\neq l\neq i$ and no sum by $i$);

\item $S_{i(j)(i)}^{l(1)(1)}=\dfrac{1}{9}\dfrac{y_{1}^{l}}{\left(
y_{1}^{i}\right) ^{2}y_{1}^{j}}$ ($i\neq j\neq l\neq i$ and no sum by $i$);

\item $S_{i(j)(k)}^{i(1)(1)}=0$ ($i\neq j\neq k\neq i$ and no sum by $i$);

\item $S_{i(l)(k)}^{l(1)(1)}=\dfrac{1}{9y_{1}^{i}y_{1}^{k}}$ ($i\neq k\neq
l\neq i$ and no sum by $l$);

\item $S_{i(j)(l)}^{l(1)(1)}=-\dfrac{1}{9y_{1}^{i}y_{1}^{j}}$ ($i\neq j\neq
l\neq i$ and no sum by $l$);

\item $S_{i(i)(l)}^{l(1)(1)}=\dfrac{1}{9\left( y_{1}^{i}\right) ^{2}}$ ($%
i\neq l$ and no sum by $i$ or $l$);

\item $S_{i(l)(i)}^{l(1)(1)}=-\dfrac{1}{9\left( y_{1}^{i}\right) ^{2}}$ ($%
i\neq l$ and no sum by $i$ or $l$);

\item $S_{l(l)(k)}^{l(1)(1)}=0$ ($k\neq l$ and no sum by $l$);

\item $S_{l(j)(l)}^{l(1)(1)}=0$ ($j\neq l$ and no sum by $l$).
\end{enumerate}
\end{theorem}

\begin{proof}
For $j\neq k$, the expression of the curvature tensor $S_{i(j)(k)}^{l(1)(1)}$
takes the form (no sum by $i$, $j$, $k$ or $l$, but with sum by $m$) 
\begin{eqnarray*}
S_{i(j)(k)}^{l(1)(1)} &=&\left[ \frac{A_{ij}^{l}\delta _{k}^{l}}{%
y_{1}^{i}y_{1}^{j}}-\frac{A_{ik}^{l}\delta _{j}^{l}}{y_{1}^{i}y_{1}^{k}}%
\right] +\left[ \frac{A_{ik}^{l}\delta _{ij}y_{1}^{l}}{\left(
y_{1}^{i}\right) ^{2}y_{1}^{k}}-\frac{A_{ij}^{l}\delta _{ik}y_{1}^{l}}{%
\left( y_{1}^{i}\right) ^{2}y_{1}^{j}}\right] + \\
&&+\left[ A_{ij}^{m}A_{mk}^{l}-A_{ik}^{m}A_{mj}^{l}\right] \frac{y_{1}^{l}}{%
y_{1}^{i}y_{1}^{j}y_{1}^{k}},
\end{eqnarray*}%
where the coefficients $A_{ij}^{l}$ are given by the relations (\ref{A-(ijk)}%
).
\end{proof}

\section{Geometrical field theories on the $1$-jet three-dimensional time}

\subsection{Geometrical gravitational theory}

\hspace{5mm}From a physical point of view, on the $1$-jet three-dimensional
time, the rheonomic Berwald-Mo\'{o}r metric (\ref{rheon-B-M}) produces the
adapted metrical d-tensor%
\begin{equation}
\mathbb{G}=h_{11}dt\otimes dt+g_{ij}dx^{i}\otimes dx^{j}+h^{11}g_{ij}\delta
y_{1}^{i}\otimes \delta y_{1}^{j},  \label{gravit-pot-B-M}
\end{equation}%
where $g_{ij}$ is given by (\ref{g-jos-(ij)}). This may be regarded as a 
\textit{"non-isotropic gravitational potential"} \cite{Mir-An}. In such a
physical context, the nonlinear connection $\mathring{\Gamma}$ (used in the
construction of the distinguished 1-forms $\delta y_{1}^{i}$) prescribes,
probably, a kind of \textit{\textquotedblleft interaction\textquotedblright }
between $(t)$-, $(x)$- and $(y)$-fields.

We postulate that the non-isotropic gravitational potential $\mathbb{G}$ is
governed by the \textit{geometrical Einstein equations}%
\begin{equation}
\text{Ric }\left( C\mathring{\Gamma}\right) -\frac{\text{Sc }\left( C%
\mathring{\Gamma}\right) }{2}\mathbb{G=}\mathcal{KT},
\label{Einstein-eq-global}
\end{equation}%
where Ric $\left( C\mathring{\Gamma}\right) $ is the \textit{Ricci d-tensor}
associated to the Cartan canonical connection $C\mathring{\Gamma}$ (in
Riemannian sense and described in adapted bases), Sc $\left( C\mathring{%
\Gamma}\right) $ is the \textit{scalar curvature}, $\mathcal{K}$ is the 
\textit{Einstein constant} and $\mathcal{T}$ is the intrinsic \textit{%
stress-energy} d-tensor of matter.

In this way, working with the adapted basis of vector fields (\ref{a-b-v}),
we can find the local geometrical Einstein equations for the rheonomic
Berwald-Mo\'{o}r metric (\ref{rheon-B-M}). Firstly, by direct computations,
we find:

\begin{proposition}
The Ricci d-tensor of the Cartan canonical connection $C\mathring{\Gamma}$
of the rheonomic Berwald-Mo\'{o}r metric (\ref{rheon-B-M}) has the following
effective adapted local Ricci d-tensors:%
\begin{equation}
\begin{array}{l}
\medskip R_{ij}=R_{ijm}^{m}=\dfrac{\varkappa _{11}^{1}\varkappa _{11}^{1}}{4}%
S_{(i)(j)}^{(1)(1)},\quad P_{i(j)}^{\text{ }(1)}=P_{(i)j}^{(1)}=P_{ij(m)}^{m%
\text{ }(1)}=\dfrac{\varkappa _{11}^{1}}{2}S_{(i)(j)}^{(1)(1)}, \\ 
S_{(i)(j)}^{(1)(1)}=S_{i(j)(m)}^{m(1)(1)}=\dfrac{3\delta _{ij}-1}{9}\cdot 
\dfrac{1}{y_{1}^{i}y_{1}^{j}}\text{ (no sum by }i\text{ or }j\text{){.}}%
\end{array}
\label{Ricci-local}
\end{equation}
\end{proposition}

\begin{remark}
The local Ricci d-tensor $S_{(i)(j)}^{(1)(1)}$ has the following expression:%
\begin{equation*}
S_{(i)(j)}^{(1)(1)}=\left\{ 
\begin{array}{ll}
-\dfrac{1}{9}\dfrac{1}{y_{1}^{i}y_{1}^{j}}, & i\neq j\medskip \\ 
\dfrac{2}{9}\dfrac{1}{\left( y_{1}^{i}\right) ^{2}}, & i=j.%
\end{array}%
\right.
\end{equation*}
\end{remark}

\begin{remark}
Using the last equality of (\ref{Ricci-local}) and the relation (\ref%
{g-sus-(jk)}), we deduce that the following equality is true (sum by $r$):%
\begin{equation}
S_{i}^{m11}\overset{def}{=}g^{mr}S_{(r)(i)}^{(1)(1)}=G_{111}^{-2/3}\cdot 
\frac{1-3\delta _{i}^{m}}{3}\cdot \frac{y_{1}^{m}}{y_{1}^{i}}\text{ (no sum
by }i\text{ or }m\text{).}  \label{S-ridicat}
\end{equation}%
Moreover, by a direct calculation, we obtain the equalities%
\begin{equation}
\sum_{m,r=1}^{3}S_{r}^{m11}C_{i(m)}^{r(1)}=0,\quad \sum_{m=1}^{3}\frac{%
\partial S_{i}^{m11}}{\partial y_{1}^{m}}=\frac{2}{3}\cdot \dfrac{1}{%
y_{1}^{i}}\cdot G_{111}^{-2/3}.  \label{equalities-S-ridicat}
\end{equation}
\end{remark}

\begin{proposition}
The scalar curvature of the Cartan canonical connection $C\mathring{\Gamma}$
of the rheonomic Berwald-Mo\'{o}r metric (\ref{rheon-B-M}) is given by%
\begin{equation*}
\text{Sc }\left( C\mathring{\Gamma}\right) =-\frac{4h_{11}+\varkappa
_{11}^{1}\varkappa _{11}^{1}}{2}\cdot G_{111}^{-2/3}.
\end{equation*}
\end{proposition}

\begin{proof}
The general formula for the scalar curvature of a Cartan connection is (for
more details, please see \cite{Neagu-Rheon})%
\begin{equation*}
\text{Sc }\left( C\mathring{\Gamma}\right)
=g^{pq}R_{pq}+h_{11}g^{pq}S_{(p)(q)}^{(1)(1)}.
\end{equation*}
\end{proof}

Describing the global geometrical Einstein equations (\ref%
{Einstein-eq-global}) in the adapted basis of vector fields (\ref{a-b-v}),
we find the following important geometrical and physical result (for more
details, please see \cite{Neagu-Rheon}):

\begin{theorem}
The adapted local \textbf{geometrical Einstein equations} that go\-vern the
non-isotropic gravitational potential (\ref{gravit-pot-B-M}), produced by
the rheonomic Berwald-Mo\'{o}r metric (\ref{rheon-B-M}), are given by:%
\begin{equation}
\left\{ 
\begin{array}{l}
\medskip \xi _{11}\cdot G_{111}^{-2/3}\cdot h_{11}=\mathcal{T}_{11} \\ 
\medskip \dfrac{\varkappa _{11}^{1}\varkappa _{11}^{1}}{4\mathcal{K}}%
S_{(i)(j)}^{(1)(1)}+\xi _{11}\cdot G_{111}^{-2/3}\cdot g_{ij}=\mathcal{T}%
_{ij} \\ 
\dfrac{1}{\mathcal{K}}S_{(i)(j)}^{(1)(1)}+\xi _{11}\cdot G_{111}^{-2/3}\cdot
h^{11}\cdot g_{ij}=\mathcal{T}_{(i)(j)}^{(1)(1)}%
\end{array}%
\right.  \label{E-1}
\end{equation}%
\medskip 
\begin{equation}
\left\{ 
\begin{array}{lll}
0=\mathcal{T}_{1i}, & 0=\mathcal{T}_{i1}, & 0=\mathcal{T}_{(i)1}^{(1)},%
\medskip \\ 
0=\mathcal{T}_{1(i)}^{\text{ }(1)}, & \dfrac{\varkappa _{11}^{1}}{2\mathcal{K%
}}S_{(i)(j)}^{(1)(1)}=\mathcal{T}_{i(j)}^{\text{ }(1)}, & \dfrac{\varkappa
_{11}^{1}}{2\mathcal{K}}S_{(i)(j)}^{(1)(1)}=\mathcal{T}_{(i)j}^{(1)},%
\end{array}%
\right.  \label{E-2}
\end{equation}%
\medskip where 
\begin{equation}
\xi _{11}=\frac{4h_{11}+\varkappa _{11}^{1}\varkappa _{11}^{1}}{4\mathcal{K}}%
.  \label{CSI}
\end{equation}
\end{theorem}

\begin{remark}
The adapted local geometrical Einstein equations (\ref{E-1}) and (\ref{E-2})
impose as the stress-energy d-tensor of matter $\mathcal{T}$ to be
symmetrical. In other words, the stress-energy d-tensor of matter $\mathcal{T%
}$ must verify the local symmetry conditions%
\begin{equation*}
\mathcal{T}_{AB}=\mathcal{T}_{BA},\quad \forall \text{ }A,B\in \left\{ 1,%
\text{ }i,\text{ }_{(i)}^{(1)}\right\} .
\end{equation*}
\end{remark}

By direct computations, the adapted local geometrical Einstein equations (%
\ref{E-1}) and (\ref{E-2}) imply the following identities of the
stress-energy d-tensor (summation by $r$):\bigskip

$\bigskip \mathcal{T}_{1}^{1}\overset{def}{=}h^{11}\mathcal{T}_{11}=\xi
_{11}\cdot G_{111}^{-2/3},\quad \mathcal{T}_{1}^{m}\overset{def}{=}g^{mr}%
\mathcal{T}_{r1}=0,\quad $

$\bigskip \mathcal{T}_{(1)1}^{(m)}\overset{def}{=}h_{11}g^{mr}\mathcal{T}%
_{(r)1}^{(1)}=0,\quad\mathcal{T}_{i}^{1}\overset{def}{=}h^{11}\mathcal{T}%
_{1i}=0,$

$\bigskip \mathcal{T}_{i}^{m}\overset{def}{=}g^{mr}\mathcal{T}_{ri}=\dfrac{%
\varkappa _{11}^{1}\varkappa _{11}^{1}}{4\mathcal{K}}S_{i}^{m11}+\xi
_{11}\cdot G_{111}^{-2/3}\cdot \mathcal{\delta }_{i}^{m},$

$\bigskip \mathcal{T}_{(1)i}^{(m)}\overset{def}{=}h_{11}g^{mr}\mathcal{T}%
_{(r)i}^{(1)}=\dfrac{h_{11}\varkappa _{11}^{1}}{2\mathcal{K}}%
S_{i}^{m11},\quad \mathcal{T}_{\text{ \ }(i)}^{1(1)}\overset{def}{=}h^{11}%
\mathcal{T}_{1(i)}^{\text{ }(1)}=0,$

$\bigskip \mathcal{T}_{\text{ \ }(i)}^{m(1)}\overset{def}{=}g^{mr}\mathcal{T}%
_{r(i)}^{\text{ }(1)}=\dfrac{\varkappa _{11}^{1}}{2\mathcal{K}}S_{i}^{m11},$

$\bigskip \mathcal{T}_{(1)(i)}^{(m)(1)}\overset{def}{=}h_{11}g^{mr}\mathcal{T%
}_{(r)(i)}^{(1)(1)}=\dfrac{h_{11}}{\mathcal{K}}S_{i}^{m11}+\xi _{11}\cdot
G_{111}^{-2/3}\cdot \mathcal{\delta }_{i}^{m},$ where the distinguished
tensor $S_{i}^{m11}$ is given by (\ref{S-ridicat}) and $\xi _{11}$ is given
by (\ref{CSI}).

\begin{theorem}
The stress-energy d-tensor of matter $\mathcal{T}$ must verify the following 
\textbf{geometrical conservation laws} (summation by $m$):%
\begin{equation*}
\left\{ 
\begin{array}{l}
\bigskip \mathcal{T}_{1/1}^{1}+\mathcal{T}_{1|m}^{m}+\mathcal{T}%
_{(1)1}^{(m)}|_{(m)}^{(1)}=\dfrac{\left( h^{11}\right) ^{2}}{16\mathcal{K}}%
\dfrac{dh_{11}}{dt}\left[ 2\dfrac{d^{2}h_{11}}{dt^{2}}-\dfrac{3}{h_{11}}%
\left( \dfrac{dh_{11}}{dt}\right) ^{2}\right] \cdot G_{111}^{-2/3} \\ 
\bigskip \mathcal{T}_{i/1}^{1}+\mathcal{T}_{i|m}^{m}+\mathcal{T}%
_{(1)i}^{(m)}|_{(m)}^{(1)}=0 \\ 
\mathcal{T}_{\text{ \ }(i)/1}^{1(1)}+\mathcal{T}_{\text{ \ }(i)|m}^{m(1)}+%
\mathcal{T}_{(1)(i)}^{(m)(1)}|_{(m)}^{(1)}=0,%
\end{array}%
\right.
\end{equation*}%
where (summation by $m$ and $r$)\bigskip

$\bigskip \mathcal{T}_{1/1}^{1}\overset{def}{=}\dfrac{\delta \mathcal{T}%
_{1}^{1}}{\delta t}+\mathcal{T}_{1}^{1}\varkappa _{11}^{1}-\mathcal{T}%
_{1}^{1}\varkappa _{11}^{1}=\dfrac{\delta \mathcal{T}_{1}^{1}}{\delta t},$

$\bigskip \mathcal{T}_{1|m}^{m}\overset{def}{=}\dfrac{\delta \mathcal{T}%
_{1}^{m}}{\delta x^{m}}+\mathcal{T}_{1}^{r}L_{rm}^{m}=\dfrac{\delta \mathcal{%
T}_{1}^{m}}{\delta x^{m}},$

$\bigskip \mathcal{T}_{(1)1}^{(m)}|_{(m)}^{(1)}\overset{def}{=}\dfrac{%
\partial \mathcal{T}_{(1)1}^{(m)}}{\partial y_{1}^{m}}+\mathcal{T}%
_{(1)1}^{(r)}C_{r(m)}^{m(1)}=\dfrac{\partial \mathcal{T}_{(1)1}^{(m)}}{%
\partial y_{1}^{m}},$

$\bigskip \mathcal{T}_{i/1}^{1}\overset{def}{=}\dfrac{\delta \mathcal{T}%
_{i}^{1}}{\delta t}+\mathcal{T}_{i}^{1}\varkappa _{11}^{1}-\mathcal{T}%
_{r}^{1}G_{i1}^{r}=\dfrac{\delta \mathcal{T}_{i}^{1}}{\delta t}+\mathcal{T}%
_{i}^{1}\varkappa _{11}^{1},$

$\bigskip \mathcal{T}_{i|m}^{m}\overset{def}{=}\dfrac{\delta \mathcal{T}%
_{i}^{m}}{\delta x^{m}}+\mathcal{T}_{i}^{r}L_{rm}^{m}-\mathcal{T}%
_{r}^{m}L_{im}^{r}=\dfrac{\varkappa _{11}^{1}}{2}\dfrac{\partial \mathcal{T}%
_{i}^{m}}{\partial y_{1}^{m}},$

$\bigskip \mathcal{T}_{(1)i}^{(m)}|_{(m)}^{(1)}\overset{def}{=}\dfrac{%
\partial \mathcal{T}_{(1)i}^{(m)}}{\partial y_{1}^{m}}+\mathcal{T}%
_{(1)i}^{(r)}C_{r(m)}^{m(1)}-\mathcal{T}_{(1)r}^{(m)}C_{i(m)}^{r(1)}=\dfrac{%
\partial \mathcal{T}_{(1)i}^{(m)}}{\partial y_{1}^{m}},$

$\bigskip \mathcal{T}_{\text{ \ }(i)/1}^{1(1)}\overset{def}{=}\dfrac{\delta 
\mathcal{T}_{\text{ \ }(i)}^{1(1)}}{\delta t}+2\mathcal{T}_{\text{ \ }%
(i)}^{1(1)}\varkappa _{11}^{1},$

$\bigskip \mathcal{T}_{\text{ \ }(i)|m}^{m(1)}\overset{def}{=}\dfrac{\delta 
\mathcal{T}_{\text{ \ }(i)}^{m(1)}}{\delta x^{m}}+\mathcal{T}_{\text{ \ }%
(i)}^{r(1)}L_{rm}^{m}-\mathcal{T}_{\text{ \ }(r)}^{m(1)}L_{im}^{r}=\dfrac{%
\varkappa _{11}^{1}}{2}\dfrac{\partial \mathcal{T}_{\text{ \ }(i)}^{m(1)}}{%
\partial y_{1}^{m}},$

$\mathcal{T}_{(1)(i)}^{(m)(1)}|_{(m)}^{(1)}\overset{def}{=}\dfrac{\partial 
\mathcal{T}_{(1)(i)}^{(m)(1)}}{\partial y_{1}^{m}}+\mathcal{T}%
_{(1)(i)}^{(r)(1)}C_{r(m)}^{m(1)}-\mathcal{T}%
_{(1)(r)}^{(m)(1)}C_{i(m)}^{r(1)}=\dfrac{\partial \mathcal{T}%
_{(1)(i)}^{(m)(1)}}{\partial y_{1}^{m}}.$
\end{theorem}

\begin{proof}
The conservation laws are provided by direct computations, using the
relations (\ref{equalitie-C}) and (\ref{equalities-S-ridicat}).
\end{proof}

\subsection{Geometrical electromagnetic theory}

\hspace{5mm}In the paper \cite{Neagu-Rheon}, using only a given Lagrangian
function $L(t,x,y)$ on the 1-jet space $J^{1}(\mathbb{R},M^{n})$, a
geometrical theory for electromagnetism was also created. In the background
of our geometrical electromagnetism from \cite{Neagu-Rheon}, we work with an 
\textit{electromagnetic distinguished }$2$\textit{-form} (the latin letters
run from $1$ to $n$)%
\begin{equation*}
\mathbb{F}=F_{(i)j}^{(1)}\delta y_{1}^{i}\wedge dx^{j},
\end{equation*}%
where%
\begin{equation*}
F_{(i)j}^{(1)}=\frac{h^{11}}{2}\left[
g_{jm}N_{(1)i}^{(m)}-g_{im}N_{(1)j}^{(m)}+\left(
g_{ir}L_{jm}^{r}-g_{jr}L_{im}^{r}\right) y_{1}^{m}\right] ,
\end{equation*}%
which is characterized by the following \textit{geometrical Maxwell equations%
} \cite{Neagu-Rheon}:%
\begin{eqnarray*}
F_{(i)j/1}^{(1)} &=&\frac{1}{2}\mathcal{A}_{\left\{ i,j\right\} }\left\{ 
\overline{D}%
_{(i)1|j}^{(1)}-D_{(i)m}^{(1)}G_{j1}^{m}+d_{(i)(m)}^{(1)(1)}R_{(1)1j}^{(m)}-%
\right. \\
&&\left. -\left[ C_{j(m)}^{p(1)}R_{(1)1i}^{(m)}-G_{i1|j}^{p}\right]
h^{11}g_{pq}y_{1}^{q}\right\} ,
\end{eqnarray*}%
\begin{equation*}
\sum_{\{i,j,k\}}F_{(i)j|k}^{(1)}=-\frac{1}{4}\sum_{\{i,j,k\}}\frac{\partial
^{3}L}{\partial y_{1}^{i}\partial y_{1}^{p}\partial y_{1}^{m}}\left[ {\dfrac{%
\delta N_{(1)j}^{(m)}}{\delta x^{k}}}-{\dfrac{\delta N_{(1)k}^{(m)}}{\delta
x^{j}}}\right] y_{1}^{p},
\end{equation*}%
\begin{equation*}
\sum_{\{i,j,k\}}F_{(i)j}^{(1)}|_{(k)}^{(1)}=0,
\end{equation*}%
where $\mathcal{A}_{\left\{ i,j\right\} }$ means an alternate sum, $%
\sum_{\{i,j,k\}}$ means a cyclic sum and we have%
\begin{equation*}
\overline{D}_{(i)1}^{(1)}=\frac{h^{11}}{2}\frac{\delta g_{im}}{\delta t}%
y_{1}^{m},\quad D_{(i)j}^{(1)}=h^{11}g_{ip}\left[
-N_{(1)j}^{(p)}+L_{jm}^{p}y_{1}^{m}\right] ,
\end{equation*}%
\begin{equation*}
d_{(i)(j)}^{(1)(1)}=h^{11}\left[ g_{ij}+g_{ip}C_{m(j)}^{p(1)}y_{1}^{m}\right]
,
\end{equation*}%
\begin{equation*}
\overline{D}_{(i)1|j}^{(1)}=\frac{\delta \overline{D}_{(i)1}^{(1)}}{\delta
x^{j}}-\overline{D}_{(m)1}^{(1)}L_{ij}^{m},\quad G_{i1|j}^{k}=\frac{\delta
G_{i1}^{k}}{\delta x^{j}}+G_{i1}^{m}L_{mj}^{k}-G_{m1}^{k}L_{ij}^{m},
\end{equation*}%
\begin{equation*}
F_{(i)j/1}^{(1)}=\frac{\delta F_{(i)j}^{(1)}}{\delta t}+F_{(i)j}^{(1)}%
\varkappa _{11}^{1}-F_{(m)j}^{(1)}G_{i1}^{m}-F_{(i)m}^{(1)}G_{j1}^{m},
\end{equation*}%
\begin{equation*}
F_{(i)j|k}^{(1)}=\frac{\delta F_{(i)j}^{(1)}}{\delta x^{k}}%
-F_{(m)j}^{(1)}L_{ik}^{m}-F_{(i)m}^{(1)}L_{jk}^{m},
\end{equation*}%
\begin{equation*}
F_{(i)j}^{(1)}|_{(k)}^{(1)}=\frac{\partial F_{(i)j}^{(1)}}{\partial y_{1}^{k}%
}-F_{(m)j}^{(1)}C_{i(k)}^{m(1)}-F_{(i)m}^{(1)}C_{j(k)}^{m(1)}.
\end{equation*}

For $n=3$, the rheonomic Berwald-Mo\'{o}r metric (\ref{rheon-B-M}) and the
nonlinear connection (\ref{nlc-B-M}), we find the electromagnetic $2$-form 
\begin{equation*}
\mathbb{F}:=\mathbb{\mathring{F}}=0.
\end{equation*}

In conclusion, our Berwald-Mo\'{o}r geometrical electromagnetic theory on
the $1$-jet three-dimensional time is trivial. In our opinion, this fact
suggests that the geometrical structure of the $1$-jet three-dimensional
time contains rather gravitational connotations than electromagnetic ones.
This is because, in our geometrical approach, the Berwald-Mo\'{o}r
electromagnetism of order three does not exist.

Gheorghe A{\scriptsize TANASIU} and Mircea N{\scriptsize EAGU}

University Transilvania of Bra\c{s}ov, Faculty of Mathematics and Informatics

Department of Algebra, Geometry and Differential Equations

B-dul Iuliu Maniu, No. 50, BV 500091, Bra\c{s}ov, Romania.

\textit{E-mails}: gh\_atanasiu@yahoo.com, mircea.neagu@unitbv.ro

\textit{Website}: http://cs.unitbv.ro/website/membrii/aged/index.html

\end{document}